\documentclass[12pt, reqno]{amsart}
\usepackage{amsfonts}
\usepackage{graphicx}
\usepackage{subfigure}
\usepackage{latexsym}
\usepackage{amsmath}
\usepackage{amssymb}
\usepackage{amsthm,amsfonts,arydshln,enumerate}
\usepackage[usenames]{color} 
\usepackage{soul} 
\usepackage{pstricks}
\usepackage{amscd}
\usepackage{color}
\numberwithin{equation}{section}

\newtheorem{thm}{Theorem}[section]
\newtheorem{lemma}[thm]{Lemma}

\newtheorem{cor}[thm]{Corollary}

\newtheorem{de}[thm]{Definition}
\newtheorem{re}[thm]{Remark}

\newtheorem*{co*}{Conjecture}
\newtheorem*{thm*}{Theorem}
\newtheorem{qu}[thm]{Question}
\def\R{\mathbb R}

\def\Z{\mathbb Z}

\begin{document}
\title{Spectrality of a class of moran measures on $\mathbb{R}^2$}

\author{Jing-Cheng Liu} \address{Key Laboratory of Computing and Stochastic Mathematics (Ministry of Education), School of Mathematics and Statistics, Hunan Normal University, Changsha, Hunan 410081, P.R. China} \email{jcliu@hunnu.edu.cn}

\author{Qiao-Qin Liu} \address{Key Laboratory of Computing and Stochastic Mathematics (Ministry of Education), School of Mathematics and Statistics, Hunan Normal University, Changsha, Hunan 410081, P.R. China} \email{lqq9991205@163.com}

\author{Jun Jason Luo} \address{College of Mathematics and Statistics, Key Laboratory of Nonlinear Analysis and its Applications (Chongqing University), Ministry of Education,  Chongqing University, Chongqing, 401331, P.R. China} \email{jun.luo@cqu.edu.cn}

\author{Jia-jie Wang} \address{Key Laboratory of Computing and Stochastic Mathematics (Ministry of Education), School of Mathematics and Statistics, Hunan Normal University, Changsha, Hunan 410081,  P.R. China} \email{wjj2021hnsd@163.com}


\thanks{The research is supported in part by the NNSF of China (No.12071125), the Hunan Provincial NSF (No.2024JJ3023), the Hunan Provincial Education Department Important Foundation of Hunan province in China (No.23A0059), the Natural Science Foundation of Chongqing (No.CSTB2023NSCQ-MSX0553).}



\begin{abstract}
	We investigate spectral properties of planar Moran measures $\mu_{\{M_n\},\{D_n\}}$ generated by sequences of expanding matrices $\{M_n\}\subset GL(2,\mathbb{Z})$ and digit sets $\{D_n\}\subset\mathbb{Z}^2$, where each digit set has the form
	\[
	D_n = \left\{
	\begin{pmatrix} 0 \\ 0 \end{pmatrix},
	\begin{pmatrix} \alpha_{n_1} \\ \alpha_{n_2} \end{pmatrix},
	\begin{pmatrix} \beta_{n_1} \\ \beta_{n_2} \end{pmatrix},
	\begin{pmatrix} -\alpha_{n_1}-\beta_{n_1} \\ -\alpha_{n_2}-\beta_{n_2} \end{pmatrix}
	\right\}
	\]
	satisfying $\alpha_{n_1}\beta_{n_2}-\alpha_{n_2}\beta_{n_1} \ne 0 \pmod{2}$.	
	Under the hypotheses $|\det(M_n)| > 4$ for all $n\geq 1$, $\sup_{n\geq 1}\|M_n^{-1}\| < 1$, and $\{D_n\}$ is finite, we establish the following characterization:
	\[
	\mu_{\{M_n\},\{D_n\}} \text{ is a spectral measure} \Longleftrightarrow M_n \in GL(2,2\mathbb{Z}) \text{ for all } n\geq 2.
	\]
	Furthermore, for the critical case $|\det(M_n)| = 4$, we derive a complete spectral criterion for  a significant class of Moran measures through combinatorial analysis of digit sets. These results extend current understanding of spectral self-affine measures to  Moran-type  constructions.
\end{abstract}

\maketitle

\section{Introduction}

A Borel measure $\mu$ supported on $\mathbb{R}^n$ is called a \emph{spectral measure} if there exists a discrete countable set $\Lambda\subset\mathbb{R}^n$ (called a \emph{spectrum}) such that the exponential system $\{e^{2\pi i\langle\lambda,x\rangle}\}_{\lambda\in\Lambda}$ forms an orthonormal basis for $L^2(\mu)$. This concept naturally extends the classical notion of spectral sets -- measurable domains $T\subset\mathbb{R}^n$ with positive and finite Lebesgue measure where $L^2(T)$ admits exponential orthonormal bases. The unit cube paradigm $([0,1]^n, \mathbb{Z}^n)$ exemplifies this fundamental relationship between geometry and spectral duality.

Fuglede's 1974 landmark work \cite{Fuglede1974} established that for Nikodym domains $T$, Segal's problem resolution coincides with $T$ being spectral. His influential spectral set conjecture proposed the equivalence:
\[
\text{Translational tiling set } \Longleftrightarrow \text{ Spectral set}.
\]
Though disproven in $\mathbb{R}^n$ ($n\geq3$) by counterexamples \cite{Kolountzakis-M2006-1,Kolountzakis-M2006-2,Matolcsi2005,Tao2004}, this conjecture remains pivotal through three key research trajectories. First, the conjecture's validity persists as open for $n=1,2$. Second, under some constrained assumptions, this conjecture has also attracted a lot of attention, e.g.,  Lev-Matolcsi's recent breakthrough \cite{Le-Ma2022} confirmed the equivalence for all convex domains. Third,  generalizing Lebesgue measure to singular measures, particularly fractal measures \cite{Gabardo-Lai2014,An-Li-Zhang2025}.

The groundbreaking work of Jorgensen-Pedersen in 1998 \cite{JP1998} marked a significant breakthrough in the field. They demonstrated that the $\frac{1}{4}$-Cantor measure is a spectral measure and explicitly identified its spectrum as
\[
\Lambda = \left\{\sum_{k=0}^n 4^k l_k : l_k\in\{0,1\},\ n\in\mathbb{N}\right\}.
\]
Strichartz \cite{Strichartz2000,Strichartz2006} subsequently expanded this to a broader class of fractal measures, including self-similar and self-affine  measures.  Laba-Wang \cite{Laba-Wang2002}, Dutkay-Jorgensen \cite{Dutkay-Jorgensen2007}, Dutkay \emph{et al.} \cite{Dutkay-Han-Sun2009}, and Dai \cite{Dai2016} characterized the spectra of spectral Cantor measures on $\R$.  Hu-Lau \cite{Hu-Lau2008}, Dai \emph{et al.} \cite{Dai2012,Dai-He-Lau2014} investigated carefully  the spectrality of Bernoulli convolutions.  More progress along this line can be found in \cite{An-He-Lau2015,An-He-Li2015,Dutkay-Lai2017,Fu-He-Wen2018,Fu-Tang2024,Li-Mi-Wa2022,Li-Mi-Wa2024,Li-Wang2024}.  On $\R^2$, the distinction between  spectral and non-spectral self-affine measures remains a central topic of ongoing research \cite{Li2008,Liu-Dong-Li2017}. In particular, a large body of work has focused on specific classes of self-affine measures, with Sierpinski-type measures receiving considerable attention \cite{Chen-Liu-Wang2024,Dai-Fu-Yan2021,Deng2016,Deng-Lau2015,Li-Zh-Wa-Ch2021,Lu-Dong-Liu2022,Wang2020}.

A matrix $M\in GL(n, \R)$ is said to be expanding if its all eigenvalues have moduli greater than one. Given an expanding matrix $M$ and a finite digit set $D\subset \R^n$, we can define an \emph{iterated function system (IFS)} $\{f_d\}_{d\in D}$ on $\mathbb{R}^n$ by
\begin{equation*}
f_{d}(x)=M^{-1}(x+d),\quad d\in D.
\end{equation*}
Then there exists a unique Borel probability measure $\mu:=\mu_{M,D}$ \cite{Hutchinson1981}  satisfying
\begin{equation}\label{eq00}
\mu_{M,D}=\frac{1}{\#D}\sum\limits_{d\in D}\mu_{M,D}\circ f_{d}^{-1}
\end{equation}
where $\#$ denotes the cardinality of a set. We call $\mu_{M,D}$ a \emph{self-affine measure}. Clearly,  $\mu_{M,D}$ is supported on the compact invariant set of the IFS.

To construct spectral self-affine measures,   Strichartz  \cite{Strichartz2000} introduced a useful concept known as \emph{admissible pairs}.

\begin{de}\label{de1.2}
Let $M\in GL(n, \Z)$ be an expanding  matrix and $D, L\subset\mathbb{Z}^n$ be two finite digit sets with the same cardinality (i.e., $\#D=\#L$). We say that ($M, D$) is an admissible pair (or ($M,D,L$) is a Hadamard triple) if the matrix
\begin{equation*}
H=\frac{1}{\sqrt{\#D}}\left(e^{2\pi i\langle M^{-1}d,\ell\rangle}\right)_{d\in D,\ell\in L}
\end{equation*}
is unitary, i.e., $H^{*}H=I$,  where $H^{*}$ means the transposed conjugate of $H$ and $I$ is the identity matrix.
\end{de}

Recently, Dutkay \emph{et al.} \cite{Dutkay-Haussermann-Lai2019} proved that if $(M,D)$ is an admissible pair, then $\mu_{M,D}$ is a spectral measure.  But the converse is usually not true, for example, $\mu_{M,D}$ is a spectral measure on $\R$ for $M=4$ and $D=\{0,1,8,9\}$, however, $(M,D)$ is not an admissible pair.   For the planar  Sierpinski-type self-affine measures,   the references \cite{Chen-Liu2019,Liu-Wang2023,Li-Zh-Wa-Ch2021} gave a complete characterization on their spectrality. More recently, Chen \emph{et al.} \cite{Chen-Liu-Wang2024} gave some necessary and sufficient conditions  for a class of self-affine measures  $\mu_{M,D}$ on $\R^2$ to be spectral measures.

\begin{thm}[\cite{Chen-Liu-Wang2024}]\label{th1.2.9}
Let $\mu_{M,D}$ be as in \eqref{eq00} where $M\in GL(2, \Z)$ is expanding and $D$ is given by
\begin{equation*}
	D=\left\{\begin{pmatrix}
		0 \\
		0
	\end{pmatrix},
	\begin{pmatrix}
		\alpha_1 \\
		\alpha_2
	\end{pmatrix},
	\begin{pmatrix}
		\beta_1 \\
		\beta_2
	\end{pmatrix},
	\begin{pmatrix}
		-\alpha_1-\beta_1 \\
		-\alpha_2-\beta_2
	\end{pmatrix}
	\right\},    \text{ where } \alpha_1\beta_2-\alpha_2\beta_1 \ne 0 \pmod{2}.
\end{equation*}
 Then the following statements are equivalent:
\begin{enumerate}
\item $\mu_{M,D}$ is a spectral measure;
\item $M\in GL(2, 2\Z)$;
 \item $(M,D)$ is an admissible pair.
\end{enumerate}
\end{thm}

Motivated by previous studies, this paper aims to extend the understanding of spectral self-affine measures (Theorem \ref{th1.2.9}) to   Moran measures, which generalize self-affine measures, possess an infinite product structure. To date, spectral analysis of Moran measures has primarily focused on the one-dimensional setting $\R$ (see \cite{An-Fu-Lai2019,An-He2014,Deng-Li2022,Li-Mi-Wa2024,Lu-Dong-Zhang2022,Luo-Mao-Liu2024,Wu-Xiao2022} and references therein). In higher dimensions, however, results concerning necessary and sufficient conditions for the spectrality of Moran measures remain limited. Notably, some progress has been made under the Hadamard triple framework \cite{Li-Wang2024,Liu-Lu-Zhou2022,Lu-Dong2021,Wang-Dong2021}, and a complete characterization was recently obtained in \cite{Deng-He-Li-Ye2024} for the spectrality of Moran-Sierpinski measures with digit set $D=\{(0,0)^t, (1,0)^t, (0,1)^t\}$.

Throughout this paper, the Moran measure under consideration is always associated with the following Moran IFS:
$$
f_{n,d}(x)=M_{n}^{-1}(x+d),\quad d\in D_n,\, n\ge 1 ,
$$
where $\left\{M_{n}\right\}\subset GL(2, \Z)$ is a sequence of expanding matrices and $\{D_n\}\subset\mathbb{Z}^{2}$ is a sequence of digit sets.  Each digit set takes the form:
\begin{equation}\label{eq1.3}
D_n=\left\{\begin{pmatrix}
0 \\
0
\end{pmatrix},\begin{pmatrix}
\alpha_{n_1} \\
\alpha_{n_2}
\end{pmatrix},\begin{pmatrix}
\beta_{n_1} \\
\beta_{n_2}
\end{pmatrix},\begin{pmatrix}
-\alpha_{n_1}-\beta_{n_1} \\
-\alpha_{n_2}-\beta_{n_2}
\end{pmatrix}
\right\},
\end{equation} where the vectors satisfy  the  determinant condition $\alpha_{n_1}\beta_{n_2}-\alpha_{n_2}\beta_{n_1} \ne 0 \pmod{2}$.

A discrete measure supported on a finite set $E\subset \R^2$ is defined by $$\delta_E=\frac{1}{\#E}\sum\limits_{e\in E}\delta_{e}$$ where $\delta_e$ is the point mass at $e$. For each integer $n\ge 1$, we write the finite convolution of discrete measures as
\begin{equation*}
\mu_{n}=\delta_{M_{1}^{-1}D_1}*\delta_{M_{1}^{-1}M_{2}^{-1}D_2}*\cdots*\delta_{M_{1}^{-1}M_{2}^{-1}\cdots M_{n}^{-1}D_n}
\end{equation*}
where the notation $*$ stands for the convolution of measures. If the sequence $\{\mu_n\}$ converges weakly to a Borel probability measure, then the weak limit is called the infinite convolution product measure, and denoted by
\begin{equation}\label{eq1.2}	\mu_{\{M_n\},\{D_n\}}=\delta_{M_{1}^{-1}D_1}*\delta_{M_{1}^{-1}M_{2}^{-1}D_2}*\cdots*\delta_{M_{1}^{-1}M_{2}^{-1}\cdots M_{n}^{-1}D_n}*\cdots.
\end{equation}
 Such  $\mu_{\{M_n\},\{D_n\}}$ (if exists)  is usually called a \emph{Moran measure}.

 It is known  \cite{Li-Mi-Wa2022}  that if
\begin{equation}\label{eq1.1.1}
\sum_{n=1}^{\infty}\max_{d\in D_n}\{\Vert M_1^{-1}M_2^{-1}\cdots M_n^{-1}d\Vert_2\}<\infty
\end{equation}
 where $\Vert \cdot \Vert_2$ is the Euclidean norm, then  the Moran measure  $\mu_{\{M_n\},\{D_n\}}$ always exists.

We define the operator norm of a matrix $A \in GL(2, \mathbb{Z})$ by
$$
\Vert A \Vert = \sup_{x \neq \mathbf{0}} \frac{\Vert Ax \Vert_2}{\Vert x \Vert_2},
$$
where $\Vert \cdot \Vert_2$ denotes the Euclidean norm. It follows that
$$
\Vert Ax \Vert_2 \leq \Vert A \Vert \cdot \Vert x \Vert_2.
$$

One of the main results of this paper is to establish a necessary and sufficient condition under which the Moran measure $\mu_{\{M_n\}, \{D_n\}}$ is spectral.

\begin{thm}\label{th1.4}
	Let $\{M_n\} \subset GL(2, \mathbb{Z})$ be a sequence of expanding matrices and $\{D_n\} \subset \mathbb{Z}^2$ a sequence of digit sets defined by \eqref{eq1.3}.  Suppose the following conditions are satisfied:
\begin{enumerate}[(i)]
\item $|\det(M_n)| > 4$ for all $n \geq 1$;

\item $\sup_{n \geq 1} \Vert M_n^{-1} \Vert < 1$;

\item   $\{D_n : n \geq 1\}$ is finite.
\end{enumerate}
Then the measure $\mu_{\{M_n\}, \{D_n\}}$ defined in \eqref{eq1.2} is a spectral measure if and only if $M_n \in GL(2, 2\mathbb{Z})$ for all $n \geq 2$.
\end{thm}

The sufficiency part of Theorem \ref{th1.4} follows directly from \cite[Theorem 1.5]{Liu-Lu-Zhou2022}. The necessity, however, is more subtle and will be derived by establishing the following more general result.

\begin{thm}\label{th1.1}
	Let $\{M_n\} \subset GL(2, \Z)$ be a sequence of expanding matrices and $\{D_n\} \subset \Z^2$  a sequence of digit sets defined by \eqref{eq1.3}. Suppose that $|\det(M_n)| \geq 4$ for all $n \geq 1$, and that the sequence $\{p_n\}$ is bounded, where $p_n = \alpha_{n_1} \beta_{n_2} - \alpha_{n_2} \beta_{n_1}$. If the Moran measure $\mu_{\{M_n\}, \{D_n\}}$ defined in \eqref{eq1.2} is a spectral measure, then it must hold that $M_n \in GL(2, 2\Z)$ for all $n \geq 2$.
\end{thm}

In the special case where $|\det(M_n)| = 4$ for all $n \geq 1$, we also establish a necessary and sufficient condition for the spectrality of a significant class of Moran measures through combinatorial analysis of digit sets. The proof is inspired by the work of Wu and Xiao \cite{Wu-Xiao2022}.

\begin{thm}\label{th1.5}
	Let $\{M_n\} \subset GL(2, 2\Z)$ be a sequence of expanding matrices with $|\det(M_n)| = 4$, and let $\Omega = \{1, 2, \dots, m\}^{\mathbb{N}}$ for some $m \geq 2$. Assume that each digit set is given by $D_n = t_n \mathcal{D}$, where $t_n \in 2\mathbb{Z} + 1$, $\mathcal{D} = \{(0,0)^t, (1,0)^t, (0,1)^t, (-1,-1)^t\}$, and $1 = t_1 < t_2 < \dots < t_m$ are pairwise coprime integers. Then, for any sequence $\sigma = (\sigma_n)_{n=1}^{\infty} \in \Omega$, the infinite convolution measure	
	$$
	\mu_\sigma = \delta_{M_1^{-1} D_{\sigma_1}} * \delta_{M_1^{-1} M_2^{-1} D_{\sigma_2}} * \cdots
	$$	
	is a spectral measure if and only if $\sigma \notin \bigcup_{l=1}^{\infty} \Sigma_l$, where $\Sigma_l = \{i_1 i_2 \cdots i_l j^\infty \in \Omega : i_l \neq j,\ j \neq 1\}$.
\end{thm}

It is worth noting that if $\sigma \in \Sigma_l$ in Theorem \ref{th1.5}, then the infinite convolution $\mu_{\sigma}$ cannot be a spectral measure. This follows from the result below.

\begin{thm}\label{th1.6}
	Let $M_1 \in GL(2, \Z)$ and $M_2 \in GL(2, 2\Z)$ be expanding matrices with $|\det(M_2)| = 4$, and define digit sets by	
	$$
	D_n = \begin{cases}
		t_1 \mathcal{D} & \text{if } n = 1, \\
		t_2 \mathcal{D} & \text{if } n \geq 2,
	\end{cases}
	\quad \text{where } t_1, t_2 \in 2\mathbb{Z} + 1,
	$$	
	and $\mathcal{D} = \{(0,0)^t, (1,0)^t, (0,1)^t, (-1,-1)^t\}$. Assume $M_n \equiv M_2$ for all $n \geq 2$. Then the Moran measure $\mu_{\{M_n\},\{D_n\}}$ defined in \eqref{eq1.2} is a spectral measure if and only if $t_2 \mid t_1$.
\end{thm}

The remainder of the paper is organized as follows. In Section \ref{sect.2}, we introduce fundamental concepts and preliminary results concerning the Fourier transform of Moran measures. Section \ref{sect.3} is devoted to the proofs of Theorems \ref{th1.4} and \ref{th1.1}. In Section \ref{sect.4}, we prove Theorems \ref{th1.5} and \ref{th1.6}. Finally, in Section \ref{sect.5}, we provide some remarks and pose open questions for future investigation.

\bigskip

\section{Preliminaries\label{sect.2}}
In this section, we will introduce some fundamental concepts and properties related to spectral measures. For a Borel probability measure $\mu$ on $\mathbb{R}^n$, the Fourier transform of $\mu$ is defined by
$$\widehat{\mu}(\xi)=\int_{\mathbb{R}^n}e^{2\pi i\langle x,\xi\rangle}\,\textup{d}\mu(x),\quad \xi\in\mathbb{R}^n.$$

By using the iteration rule, the Fourier transform of  the self-affine measure $\mu_{M,D}$ in \eqref{eq00} can be written as
\begin{equation}\label{eq2.1}
\widehat{\mu}_{M,D}(\xi)=\prod_{j=1}^\infty m_D({M^{*}}^{-{j}}\xi)
\end{equation}
where $m_{D}(\cdot)=\frac{1}{\#D}\sum\limits_{d\in D}e^{2\pi i\langle d,\cdot\rangle}$, which is named by the mask polynomial of the digit set $D$. It is obvious that $m_{D}(\cdot)$ is a $\mathbb{Z}^n$-periodic function when $D\subset\mathbb{Z}^n$.

Similarly, the Fourier transform of the Moran measure $\mu_{\{M_n\},\{D_n\}}$ in \eqref{eq1.2} is of the form
\begin{equation}\label{eq2.2}
\widehat{\mu}_{\{M_n\},\{D_n\}}(\xi)=\prod\limits_{j=1}^{\infty}m_{D_j}((M_{1}^{*}M_{2}^{*}\cdots M_{j}^{*})^{-1}\xi).
\end{equation}

If we let $\mathcal{Z}(f)=\{x:f(x)=0\}$ denote the zero set of a function $f$, then it directly follows from \eqref{eq2.1} and \eqref{eq2.2} that
\begin{equation*}
\mathcal{Z}(\widehat{\mu}_{M,D})=\bigcup_{j=1}^{\infty}M^{*j}\mathcal{Z}(m_{D})
\end{equation*}
and
\begin{equation}\label{eq2.4}
\mathcal{Z}(\widehat{\mu}_{\{M_n\},\{D_n\}})=\bigcup_{j=1}^{\infty}M_{1}^{*}M_{2}^{*}\cdots M_{j}^{*}\mathcal{Z}(m_{D_j}).
\end{equation}

Let $\Lambda$ be a countable subset of $\R^n$, then $E_\Lambda:=\{e^{2\pi i\langle\lambda,x\rangle}:\lambda\in\Lambda\}$ forms an orthogonal family for $L^2(\mu)$ if and only if $\widehat{\mu}(\lambda_1-\lambda_2)=0$ holds for any $\lambda_1\neq\lambda_2\in E_{\Lambda}$, i.e.,
\begin{equation}\label{eq2.5}
 (\Lambda-\Lambda)\setminus\{0\}\subset\mathcal{Z}(\widehat{\mu}).
\end{equation}
We call $\Lambda$ an \emph{orthogonal set (respectively, spectrum)} of $\mu$ if $E_\Lambda$ forms an orthogonal family (respectively, Fourier basis) for $L^2(\mu)$. The following lemma indicates that the spectrality of   Moran measure in \eqref{eq1.2} is invariant under a similarity transformation. We omit the proof here as it is analogous to \cite[Lemma 4.2]{Dutkay-Jorgensen2007}.

\begin{lemma}\label{lem2.2}
Let $\{D_n\}, \{\widetilde{D}_n\}$ be two sequences of finite digit sets and $\{M_n\}, \{\widetilde{M}_n\}$ be two sequences of expanding matrices. If there exists an invertible matrix  $Q$ such that $\widetilde{M}_n=QM_nQ^{-1}$ and $\widetilde{D}_n=QD_n$ for all $n\geq 1$, then $\mu_{\{M_n\},\{D_n\}}$ is a spectral measure with spectrum $\Lambda$ if and only if $\mu_{\{\widetilde{M}_n\},\{\widetilde{D}_n\}}$ is a spectral measure with spectrum $Q^{*-1}\Lambda$.
\end{lemma}

The next theorem is a basic criterion for an orthogonal set to be a spectrum of a measure. Denote by $$Q_{\mu,\Lambda}(\xi):=\sum\limits_{\lambda\in\Lambda}|\widehat{\mu}(\xi+\lambda)|^2, \quad \xi\in \R^n.$$

\begin{thm}[\cite{JP1998}]\label{thm-JP98}
	Let $\mu$ be a Borel probability measure with compact support on $\mathbb{R}^n$, and $\Lambda\subset\mathbb{R}^n$ a countable set. Then
\begin{enumerate}
\item $\Lambda$ is an orthogonal set of $\mu$ if and only if $Q_{\mu,\Lambda}(\xi)\leq1$ for all $\xi\in\mathbb{R}^n$;
\item $\Lambda$ is a spectrum of $\mu$ if and only if $Q_{\mu,\Lambda}(\xi)\equiv1$ for all $\xi\in\mathbb{R}^n$.
\end{enumerate}
\end{thm}

Given an expanding matrix $M$ and a digit set $D$, the following lemma provides a criterion for determining whether $(M,D)$ forms an admissible pair.
	
\begin{lemma}[\cite{Dutkay-Haussermann-Lai2019}]\label{lem2.3}
Let $M\in GL(n, \Z)$ be an expanding matrix, and let $D,L\subset\mathbb{Z}^n$ be two finite digit sets with the same cardinality. Then the following are equivalent:
\begin{enumerate}
 	\item $(M,D,L)$ is a Hadamard triple;
 	\item $m_D(M^{*-1}(l_1-l_2))=0$ for any $l_1\neq l_2\in L$;
    \item $(\delta_{M^{-1}D},L)$ is a spectral pair.
    \end{enumerate}
\end{lemma}

In the proof of Theorem \ref{th1.4}, we also require the following technical lemma, originally stated by Deng and Li \cite[Lemma 2.5]{Deng-Li2022}.

\begin{lemma}[\cite{Deng-Li2022}]\label{lem2.6}
Let $p_{i,j}$ be positive numbers such that $\sum\limits_{j=1}^{n}p_{i,j}=1$ and $q_{i,j}$ be nonnegative numbers such that $\sum\limits_{i=1}^{m}\mathop{max}\limits_{1\leq j\leq n}\{q_{i,j}\}\leq1$. Then $\sum\limits_{i=1}^{m}\sum\limits_{j=1}^{n}p_{i,j}q_{i,j}=1$ if and only if $q_{i,1}=\cdots=q_{i,n}$ for $1\leq i\leq m$ and $\sum\limits_{i=1}^{m}q_{i,1}=1$.
\end{lemma}

Let  $\mu$ be a Borel probability measure on $\mathbb{R}^n$,   the \emph{integer periodic zero set} of $\mu$ is given by
\begin{equation*}
\textbf{Z}(\mu)=\{\xi\in\mathbb{R}^n:\widehat{\mu}(\xi+k)=0 \ \text{for all} \ k\in\mathbb{Z}^n\}.
\end{equation*}
This concept was first proposed by Dutkay \emph{et al.} \cite{Dutkay-Haussermann-Lai2019}, which plays an important role in seeking spectral measures. For further analysis, we also need  a sequence of measures $\{\nu_n\}$  induced by  $\mu_{\{M_n\},\{D_n\}}$ in \eqref{eq1.2}. The $\nu_n$'s are defined by
\begin{equation}\label{eq1.2.0}
\nu_n :=\delta_{M_{n+1}^{-1}D_{n+1}}*\delta_{M_{n+1}^{-1}M_{n+2}^{-1}D_{n+2}}*\cdots
\end{equation}
for all $n\ge 1$.

Recently, Li and Wang \cite{Li-Wang2024} established a useful criterion for determining when the infinite convolution generated by a sequence of admissible pairs yields a spectral measure.

\begin{thm}[\cite{Li-Wang2024}]\label{th2.5}
Given a sequence of admissible pairs $\{(M_n,D_n)\}$ on $\mathbb{R}^n$, suppose that the infinite convolution $\mu_{\{M_n\},\{D_n\}}$ defined by \eqref{eq1.2} exists and satisfies
$$\lim_{n\rightarrow\infty}\|M_1^{-1}M_2^{-1}\cdots M_n^{-1}\|=0.$$
If there exists a subsequence $\{\nu_{n_j}\}$ given in \eqref{eq1.2.0} which converges weakly to $v$ with $\textup{\textbf{Z}}(v)=\emptyset$, then $\mu_{\{M_n\},\{D_n\}}$ is a spectral measure.
\end{thm}

In proving Theorem \ref{th1.5}, it suffices to verify the condition $\sup_{n\ge 1}\Vert M_n^{-1}\Vert<1$, which is stronger than the above one.

\bigskip

\section{Proofs of Theorems \ref{th1.4} and \ref{th1.1}\label{sect.3}}

First of all, we provide several technical lemmas and then use them to prove Theorem \ref{th1.1}. Secondly,  we  apply Theorem \ref{th1.1} to complete the proof of Theorem \ref{th1.4} at the end of this section.

Let $Q_0=I$, $Q_n:=\begin{pmatrix}
	\alpha_{n_1} & \beta_{n_1} \\
	\alpha_{n_2} & \beta_{n_2}
\end{pmatrix},  n\geq 1$, and let
$$\widetilde{M}_n=Q_{n}^{-1}M_nQ_{n-1}, \quad \mathcal {D}=\left\{\begin{pmatrix}
		0\\
		0
		\end{pmatrix},\begin{pmatrix}
		1\\
		0
		\end{pmatrix},\begin{pmatrix}
		0\\
		1
		\end{pmatrix},\begin{pmatrix}
		-1\\
		-1
		\end{pmatrix}\right\}.$$

Since $D_n=Q_n\mathcal {D}$ and  $M_1^{-1}\cdots M_n^{-1}D_n=\widetilde{M}_1^{-1}\cdots \widetilde{M}_n^{-1}\mathcal{D}$, we can rewrite the Moran measure $\mu_{\{M_n\},\{D_n\}}$ defined in \eqref{eq1.2} as the following form:
\begin{equation}\label{eq3.1}
\begin{split}
\mu_{\{M_n\},\{D_n\}}&= \delta_{\widetilde{M}_1^{-1}\mathcal {D}}*\delta_{\widetilde{M}_1^{-1}\widetilde{M}_2^{-1}\mathcal {D}}*\cdots*\delta_{\widetilde{M}_1^{-1}\widetilde{M}_2^{-1}\cdots \widetilde{M}_n^{-1}\mathcal {D}}*\cdots  \\
&:=\mu_{\{\widetilde{M}_n\},\mathcal {D}}.
\end{split}
\end{equation}
For simplicity, we use $\mu$ to denote $\mu_{\{\widetilde{M}_n\},\mathcal {D}}$ in the rest of this section.

According to the assumption of Theorem \ref{th1.1}, the sequence $\{p_n=\det(Q_n)\}$  is bounded. Hence $p_n$'s only take a finite number of values. Without loss of generality, we may assume that these values are $p_1,p_2,\dots,p_m$, and denote by $p=\text{lcm}(p_1,p_2,\dots,p_m)$.  For the sake of brevity, we always assume that $p\geq 2$, since the case $p=1$ can be shown in an analogous way. For $n>1$, we write
\begin{equation}\label{eq-3-2-1}
\mathcal{F} _n:=\big\{(\ell_1,\ell_2)^t:\ell_1,\ell_2\in\{0,1,\dots,n-1\}\big\} \quad\text{and}\quad \mathring{\mathcal{F}}_n:=\mathcal{F}_n\setminus\{\mathbf{0}\}.
\end{equation}

From the definitions of mask polynomials and zero sets in the preceding section, it is immediate to find the relation
\begin{equation}\label{eq3.2.0}
\mathcal{Z}(m_\mathcal{D})=\frac{1}{2}\mathring{\mathcal{F}}_2+\mathbb{Z}^2.
\end{equation}

In the following, we investigate the structure of the spectrum of the Moran measure $\mu$. Let $\Lambda$ be a spectrum of $\mu$ with $\mathbf{0}\in\Lambda$. Then by \eqref{eq2.4}, \eqref{eq2.5} and \eqref{eq3.2.0}, we have
\begin{equation}\label{eq4.6}
\Lambda\subset\{\mathbf{0}\}\cup\bigcup_{n=1}^{\infty}\widetilde{M}^{*}_{1}\widetilde{M}_{2}^{*}\widetilde{M}_{3}^{*}\cdots \widetilde{M}_{n}^{*}\mathcal{Z}(m_{\mathcal {D}}) \ \ \ {\rm and }\ \ \ 2p\widetilde{M}_{1}^{*-1}\Lambda\subset\mathbb{Z}^2.
\end{equation}

It is clear that $\mathcal{F}_p\oplus\ p\mathcal{F}_2$ is a complete residue system of the diagonal matrix $2pI$. Hence for any ${\mathbf k}\in\mathbb{Z}^2$, there exist $s\in \mathcal{F}_p, l\in \mathcal{F}_2$ and ${\mathbf k'}\in\mathbb{Z}^2$ such that
\begin{equation}\label{eq4.6.1}
{\mathbf k}=s+pl+2p{\mathbf k'}.
\end{equation}

By \eqref{eq4.6} and \eqref{eq4.6.1}, the spectrum $\Lambda$ of $\mu$ has the following decomposition:
\begin{equation}\label{eq4.7.0}
	\Lambda=\frac{\widetilde{M}_{1}^{*}}{2p}\bigcup\limits_{s\in \mathcal{F}_p}\bigcup\limits_{l\in \mathcal{F}_2}\left(s+pl+2p\Lambda_{s,l}\right)
\end{equation}
where
\begin{equation}\label{eq-3.6.0}
\Lambda_{s,l}=\left\{\gamma\in\mathbb{Z}^2:s+pl+2p\gamma\in2p\widetilde{M}_{1}^{*-1}\Lambda\right\}
\end{equation}
and $s+pl+2p\Lambda_{s,l}=\emptyset$ if $\Lambda_{s,l}=\emptyset$. In addition, $\Lambda_{\mathbf{0},\mathbf{0}}\neq\emptyset$ since $\textbf{0}\in\Lambda$.

So far, we give the structure of the spectrum $\Lambda$ of $\mu$ under the assumption that $\mu$ is a spectral measure. In the following, we turn to characterize the spectrum of the following Moran measure
\begin{equation}\label{eq4.8.1}
\mu_1=\delta_{\widetilde{M}_{2}^{-1}\mathcal {D}}*\delta_{\widetilde{M}_2^{-1}\widetilde{M}_{3}^{-1}\mathcal {D}}*\cdots.
\end{equation}

\begin{lemma}\label{lem4.2}
Let $\Lambda$ be a spectrum of $\mu$ with $\mathbf{0}\in\Lambda$. For any $s\in \mathcal{F}_p$, select an $l_{s}\in\mathcal{F}_2$ and let
\begin{equation*}
\Gamma=\bigcup\limits_{s\in \mathcal{F}_p}\left(\frac{s+pl_s}{2p}+\Lambda_{s,l_s}\right)
\end{equation*}
where $\Lambda_{s,l_s}$ is as in \eqref{eq-3.6.0}. If $\Gamma\neq\emptyset$, then $\Gamma$ is a spectrum of $\mu_1$ defined by \eqref{eq4.8.1}.
\end{lemma}

\begin{proof}
Our proof  is divided into  two steps.

\noindent\textbf{Step I.} To show that $\Gamma$ is an orthogonal set of $\mu_1$.

For any $\varsigma_1\neq\varsigma_2\in\Gamma$, it follows from the definition of $\Gamma$ that there exist $s_k\in \mathcal{F}_p$, $l_{s_k}\in\mathcal{F}_2$ and $\gamma_k\in\Lambda_{s_k,l_{s_k}}$, $k=1,2$ such that
\begin{equation*}
	\varsigma_k=\frac{s_k+pl_{s_k}}{2p}+\gamma_k.
\end{equation*}
Thus, by \eqref{eq2.5} and \eqref{eq4.7.0}, we have $\widetilde{M}_{1}^{*}(\varsigma_1-\varsigma_2)\in(\Lambda-\Lambda)\backslash\{\textbf{0}\}\subset\mathcal{Z}(\widehat{\mu})$. Noting that $\gamma_1,\gamma_2\in\mathbb{Z}^2$, the $\mathbb{Z}^2$-periodicity of $m_{\mathcal {D}}$ yields
\begin{equation*}
\begin{split}
0&=\widehat{\mu}(\widetilde{M}_{1}^{*}(\varsigma_1-\varsigma_2))=
m_{\mathcal {D}}(\varsigma_1-\varsigma_2)\widehat{\mu}_1(\varsigma_1-\varsigma_2) \\
&=m_{\mathcal {D}}\left(\frac{s_1-s_2+p(l_{s_1}-l_{s_2})}{2p}
+\gamma_1-\gamma_2\right)\widehat{\mu}_1(\varsigma_1-\varsigma_2) \\
&=m_{\mathcal {D}}\left(\frac{s_1-s_2}{2p}+\frac{l_{s_1}-l_{s_2}}{2}
\right)\widehat{\mu}_1(\varsigma_1-\varsigma_2).
\end{split}
\end{equation*}

Once we have $m_{\mathcal {D}}\left(\frac{s_1-s_2}{2p}+\frac{l_{s_1}-l_{s_2}}{2}
\right)\neq0$, we derive $\widehat{\mu}_1(\varsigma_1-\varsigma_2)=0$, which indicates that $\Gamma$ is an orthogonal set of $\mu_1$.
Next, we shall verify $m_{\mathcal {D}}\left(\frac{s_1-s_2}{2p}+\frac{l_{s_1}-l_{s_2}}{2}
\right)\neq0$ by considering the sequel two cases.

\textbf{Case i:} $s_1=s_2$.  Then  $l_{s_1}=l_{s_2}$. It is clear that
$$
m_{\mathcal {D}}\left(\frac{s_1-s_2}{2p}+\frac{l_{s_1}-l_{s_2}}{2}
\right)=m_{\mathcal {D}}(\textbf{0})=1\neq0.
$$

\textbf{Case ii:} $s_1\neq s_2$. Assume conversely that $m_{\mathcal {D}}\left(\frac{s_1-s_2}{2p}+\frac{l_{s_1}-l_{s_2}}{2}
\right)=0$. Then
\begin{equation*}
\frac{s_1-s_2}{2p}+\frac{l_{s_1}-l_{s_2}}{2}\in\mathcal{Z}(m_{\mathcal {D}})
=\frac{1}{2}(\mathbb{Z}^2\setminus 2\mathbb{Z}^2).
\end{equation*}
It yields that $s_1-s_2\in p\mathbb{Z}^2$, which is impossible as $s_1,s_2 \in \mathcal{F}_p$.

\bigskip

\noindent\textbf{Step II.} To show that $\Gamma$ is a spectrum of $\mu_1$.

Firstly, it is easy to verify that the triple ($2I,\mathcal {D},\mathcal{F}_2$) is a Hadamard triple by Lemma \ref{lem2.3} (ii).  For fixed $s\in \mathcal{F}_p$,  it follows from Lemma \ref{lem2.3} (iii) and Theorem \ref{thm-JP98} (ii) that
\begin{equation}\label{eq4.9}	
\sum\limits_{l\in \mathcal{F}_2}\left|m_{\mathcal {D}}\left(\xi+\frac{s+pl}{2p}\right)\right|^2=1
\end{equation}
holds for any $\xi\in\mathbb{R}^2$.

Because $\Lambda$ is a spectrum of $\mu$ and $\Lambda_{s,l}\subset\mathbb{Z}^2$, using Theorem \ref{thm-JP98} (ii), $\mathbb{Z}^2$-periodicity of $m_{\mathcal {D}}$, \eqref{eq4.7.0} and \eqref{eq4.9},  we can get, for $\xi\in\mathbb{R}^2\setminus\mathbb{Q}^2$,
\begin{align}\label{eq3.11}
   1&\equiv \sum_{\lambda\in\Lambda}|\widehat{\mu}(\widetilde{M}_{1}^{*} \xi+\lambda)|^2 \nonumber \\
	&=  \sum_{s\in \mathcal{F}_p}\sum_{l\in \mathcal{F}_2}\sum\limits_{\gamma\in\Lambda_{s,l}}\left|\widehat{\mu}
\left(\widetilde{M}_{1}^{*}\xi+\frac{\widetilde{M}_{1}^{*}(s+pl+2p\gamma)}{2p}\right)\right|^2\nonumber\\
	&= \sum\limits_{s\in \mathcal{F}_p}\sum\limits_{l\in \mathcal{F}_2}\left|m_{\mathcal {D}}\left(\xi+\frac{s+pl}{2p}\right)
\right|^2\sum\limits_{\gamma\in\Lambda_{s,l}}\left|\widehat{\mu}_1
\left(\xi+\frac{s+pl}{2p}+\gamma\right)\right|^2.
\end{align}
Denote
\begin{equation*} p_{s,l}=\left|m_{\mathcal {D}}\left(\xi+\frac{s+pl}{2p}\right)\right|^2
\end{equation*}
and
\begin{equation*} q_{s,l}=\sum\limits_{\gamma\in\Lambda_{s,l}}\left|\widehat{\mu}_1
\left(\xi+\frac{s+pl}{2p}+\gamma\right)\right|^2\geq0.
\end{equation*}
Clearly $p_{s,l}>0$ as $\xi\in\mathbb{R}^2\setminus\mathbb{Q}^2$ and \eqref{eq3.2.0}. Thus \eqref{eq3.11} becomes
\begin{equation}\label{eq4.11.0}
\sum\limits_{s\in \mathcal{F}_p}\sum\limits_{l\in \mathcal{F}_2}p_{s,l}q_{s,l}=1.
\end{equation}
By \eqref{eq4.9}, we have $\sum\limits_{l\in \mathcal{F}_2}p_{s,l}=1$.
Furthermore, since $\Gamma$ is an orthogonal set of $\mu_1$, it follows from Theorem \ref{thm-JP98} (i) that $\sum\limits_{s\in \mathcal{F}_p}\max\limits_{l\in\mathcal{F}_2}\{q_{s,l}\}\leq1$. Combining  \eqref{eq4.11.0} and Lemma \ref{lem2.6}, it concludes that
\begin{equation}\label{eq4.12}
\sum\limits_{s\in \mathcal{F}_p}\sum\limits_{\gamma\in\Lambda_{s,l^*}}\left|\widehat{\mu}_1
\left(\xi+\frac{s+pl^*}{2p}+\gamma\right)\right|^2=1,\ \ l^*\in\mathcal{F}_2
\end{equation}
and
\begin{eqnarray}\label{eq4.13}
\sum\limits_{\gamma\in\Lambda_{s,l_1}}\left|\widehat{\mu}_1
\left(\xi+\frac{s+pl_1}{2p}+\gamma\right)\right|^2 =\sum\limits_{\gamma\in\Lambda_{s,l_2}}\left|\widehat{\mu}_1
\left(\xi+\frac{s+pl_2}{2p}+\gamma\right)\right|^2
\end{eqnarray}
for any $s\in \mathcal{F}_p$ and $l_1,l_2\in\mathcal{F}_2$.

By the continuity of Fourier transform, \eqref{eq4.12} and \eqref{eq4.13} hold for any $\xi\in\mathbb{R}^2$. According to Theorem \ref{thm-JP98} (ii), $\Gamma$ is a spectrum of $\mu_1$ for any group $\{l_s\}_{s\in\mathcal{F}_p}$. We complete the proof of the lemma.
\end{proof}

\begin{re}\label{re3.5}
For $n\ge 1$, let
\begin{equation*}
\mu_n=\delta_{\widetilde{M}_{n+1}^{-1}\mathcal {D}}*\delta_{\widetilde{M}_{n+1}^{-1}\widetilde{M}_{n+2}^{-1}\mathcal {D}}*\cdots.
\end{equation*}

By repeatedly using Lemma \ref{lem4.2}, we actually obtain that if $\mu$ is a spectral measure, then $\mu_n$ is a spectral measure for each $n\ge 1$.
\end{re}

\bigskip

\begin{proof}[\textbf{Proof of Theorem \ref{th1.1}}]
Suppose that $\mu$ is a spectral measure with spectrum $\Lambda$ and $\mathbf{0}\in\Lambda$. Firstly, we construct a spectrum $\Gamma$ of $\mu_1$ in \eqref{eq4.8.1}.  For any $s\in\mathcal{F}_p$, we choose $l_s\in\mathcal{F}_2$ such that $s+pl_s\in 2\mathbb{Z}^2$, then Lemma \ref{lem4.2} shows that
\begin{equation*}
\Gamma=\bigcup\limits_{s\in \mathcal{F}_p}\left(\frac{s+pl_s}{2p}+\Lambda_{s,l_s}\right)
\end{equation*}
is a spectrum of $\mu_1$ and $\mathbf{0}\in\Gamma$. By \eqref{eq4.7.0}, the spectrum $\Gamma$ can also be decomposed as the following form:
\begin{equation}\label{eq4.7}
	\Gamma=\frac{\widetilde{M}_{2}^{*}}{2p}\bigcup\limits_{s'\in \mathcal{F}_p}\bigcup\limits_{l'\in \mathcal{F}_2}\left(s'+pl'+2p\Lambda'_{s',l'}\right).
\end{equation}

Using \eqref{eq4.13} to $\mu_2$, we conclude that $\Lambda'_{\mathbf{0},l'}\neq\emptyset$ for all $l'\in \mathcal{F}_2$ since $\Lambda'_{\mathbf{0}, \mathbf{0}}\neq\emptyset$. For any $l'\in \mathcal{F}_2$, choose an element
$z'_{l'}\in \Lambda'_{\mathbf{0},l'}$, then \eqref{eq4.7} implies that there exist $s\in \mathcal{F}_p$, $l_s\in\mathcal{F}_2$ and $z\in \Lambda_{s,l_s}$ such that
\begin{equation}\label{eq4.8}
s+pl_s+2pz=\widetilde{M}_{2}^{*}(pl'+2pz'_{l'})
=Q_{1}^{*}M_{2}^{*}Q_{2}^{*-1}(pl'+2pz'_{l'}).
\end{equation}
Noting that $\det(Q_2) \mid p$ and $s+pl_s\in 2\mathbb{Z}^2$, it follows from \eqref{eq4.8} that for any $l'\in \mathcal{F}_2$,
\begin{equation*}
	Q_{1}^{*}M_{2}^{*}Q_{2}^{*-1}pl'=(s+pl_s+2pz)-Q_{1}^{*}M_{2}^{*}Q_{2}^{*-1}2pz'_{l'}\in 2\mathbb{Z}^2.
\end{equation*}

Write $B_2 := pQ_{1}^{*}M_{2}^{*}Q_{2}^{*-1}$. It can be verified easily that $B_2(1,0)^t, B_2(0,1)^t, B_2(1,1)^t$ all lie in $2\mathbb{Z}^{2}$. Hence $B_2\in GL(2,2\Z)$, and $p\det(Q_{1}^{*})M_{2}^{*}=\det(Q_{1}^{*})Q_{1}^{*-1}B_2Q_{2}^{*}\in GL(2,2\Z)$.
As $p\det(Q_{1}^{*})\in 2\mathbb{Z}+1$, we prove that $M_2\in  GL(2,2\Z)$. Repeating the similar argument, we further have $M_n\in  GL(2,2\Z)$ for any $n> 2$.  That finishes the proof of Theorem \ref{th1.1}.
\end{proof}

\bigskip

\begin{proof}[\textbf{Proof of Theorem \ref{th1.4}}]
The necessity is due to Theorem \ref{th1.1}. Now we show the sufficiency. From the assumption: $\sup\limits_{n\ge 1}\Vert M_n^{-1}\Vert<1$ and $\#\{D_n:n\ge 1\}<\infty$, it follows that $\sup\limits_{n\ge 1}\sup\limits_{d\in D_n}\Vert M_n^{-1}d\Vert_2<\infty$. Moreover, as $M_n\in GL(2,2\Z)$ for any $n\geq2$,  Theorem \ref{th1.2.9} implies that $\{(M_n,D_n)\}_{n=2}^{\infty}$ forms a sequence of admissible pairs. On the other hand, by  \cite[Lemma 2.6]{Deng-He-Li-Ye2024}, the spectrality of the Moran measure $\mu_{\{M_n\},\{D_n\}}$ is independent of the choice of  $M_1$.  Therefore, we can take $M_1=4I$, which forces  $(M_1,D_1)$ to be an admissible pair as well by Theorem \ref{th1.2.9}.  Consequently,  the sufficiency follows from \cite[Theorem 1.5]{Liu-Lu-Zhou2022}.
\end{proof}

\bigskip

\section{Proofs of Theorems \ref{th1.5} and \ref{th1.6} \label{sect.4}}
The proof scheme of  Theorem \ref{th1.5} and Theorem \ref{th1.6} is as follows: Firstly, we give some basic lemmas for proving the sufficiency of Theorem \ref{th1.5} and then prove Theorem \ref{th1.6}; finally we prove  Theorem \ref{th1.5}.

Recall the assumptions of Theorem  \ref{th1.5}.  $D_n=t_n\mathcal{D}, n=1,\dots, m$ where $ t_n\in2\mathbb{Z}+1, \mathcal{D}=\{(0,0)^t,(1,0)^t,(0,1)^t,(-1,-1)^t\}$, and  $1=t_1<t_2<\dots<t_m$ are pairwise coprime numbers.  $\{M_n\}\subset GL(2, 2\Z)$ is a sequence of expanding matrices with $|\det(M_n)|=4$. It is easy to verify that the moduli of two eigenvalues of every $M_n$ are both equal to $2$. Hence  $\iota :=\sup\limits_{n\ge 1}\Vert M_n^{-1}\Vert=1/2<1$.  Write 
$$\gamma:=\max\{\Vert d\Vert_2:d\in D_n, 1\leq n\leq m\},$$ 
then $\gamma<\infty$. For any  $\sigma=(\sigma_k)_{k=1}^{\infty}\in\Omega=\{1,2,\dots,m\}^{\mathbb{N}}$, the infinite convolution  $\mu_\sigma$ is
\begin{equation}\label{eq5.1.1}
	\mu_\sigma=\delta_{M_{1}^{-1}D_{\sigma_1}}*\delta_{M_{1}^{-1}M_{2}^{-1}D_{\sigma_2}}*\cdots.
\end{equation}
Since
\begin{align}\label{eq5.4.0}
\sum_{n=1}^{\infty}\max_{d_{\sigma_n}\in D_{\sigma_n}}\{\Vert M_{1}^{-1}M_{2}^{-1}\cdots M_{n}^{-1}d_{\sigma_n}\Vert_2\}&\leq\sum_{n=1}^{\infty}\Vert M_{1}^{-1}\Vert\cdot\Vert M_{2}^{-1}\Vert\cdots\Vert M_{n}^{-1}\Vert\cdot\gamma\nonumber\\
&\leq \sum_{n=1}^{\infty}\gamma\iota^{n}=\frac{\gamma\iota}{1-\iota}=\gamma<\infty,
\end{align}
it follows from \eqref{eq1.1.1} that $\mu_{\sigma}$ exists and has a compact support denoted by $\text{spt}(\mu_{\sigma})$.

We define a metric $d$ on the symbol space $\Omega$ by $d(\sigma,\tau)=2^{-\min\{n\geq1 : \sigma_n\neq\tau_n\}}$
where $\sigma=(\sigma_n)_{n=1}^{\infty}$, $\tau=(\tau_n)_{n=1}^{\infty}\in\Omega$. Trivally, $(\Omega,d)$ is a  compact metric space. We say that a sequence $\{\sigma(j)\}_{j=1}^{\infty}\subset\Omega$ converges to $\tau\in\Omega$,   if for each $n\geq1$, there exists $j_0\geq1$ such that for all $j\geq j_0$,
$$\sigma_1(j)\sigma_2(j)\cdots \sigma_n(j)=\tau_1\tau_2\cdots\tau_n.$$

The convergence of sequences in $\Omega$ yields the weak convergence of the sequences of the corresponding measures due to the following result. The proof is analogous to \cite[Lemma 5.1]{Li-Mi-Wa2024}.

\begin{lemma}\label{lem5.2.0}
Given a sequence $\{\sigma(j)\}_{j=1}^{\infty}\subset\Omega$ and $\tau\in\Omega$,  if $\{\sigma(j)\}_{j=1}^{\infty}$ converges to $\tau$, then
$\{\mu_{\sigma(j)}\}_{j=1}^{\infty}$ defined by \eqref{eq5.1.1} converges weakly to $\mu_{\tau}$.
\end{lemma}
\begin{proof}
For $n\ge 1$ and $\sigma\in\Omega$, we define
\begin{equation*}
\mu_{\sigma,n}=\delta_{M_{1}^{-1}D_{\sigma_1}}*\delta_{M_{1}^{-1}
M_{2}^{-1}D_{\sigma_2}}*\cdots*\delta_{M_{1}^{-1}\cdots M_{n}^{-1}D_{\sigma_n}}
\end{equation*}
and
\begin{equation*}
\mu_{\sigma,>n}=\delta_{M_{1}^{-1}\cdots M_{{n+1}}^{-1}D_{\sigma_{n+1}}}*\delta_{M_{1}^{-1}\cdots M_{{n+2}}^{-1}D_{\sigma_{n+2}}}*\cdots.
\end{equation*}
It is clear that $\mu_{\sigma}=\mu_{\sigma,n}*\mu_{\sigma,>n}$. Let $B(x,r)$ denote the closed ball with center  $x\in\mathbb{R}^2$ and radius $r$. From \eqref{eq5.4.0}, it can be seen that the supports of the measures satisfy
$$\text{spt}(\mu_{\sigma,n})\subset B(\textbf{0},\gamma),\quad \text{spt}(\mu_{\sigma,>n})\subset B(\textbf{0},\frac{\gamma}{2^n}).$$

Let $f$ be any bounded continuous function on $\mathbb{R}^2$, then $f$ is uniformly continuous on the closed ball $B(\textbf{0},\frac{\gamma}{2^n})$. Thus, for any $\varepsilon>0$, there exists $\delta>0$ such that for all $x,y\in  B(\textbf{0},\frac{\gamma}{2^n})$, we have
\begin{equation*}
|f(x)-f(y)|<\varepsilon \quad\text{provided}\quad \Vert x-y \Vert_2<\delta.
\end{equation*}
Take a sufficiently large $n_{0}$ such that $\frac{\gamma}{2^{n_0}}<\frac{\delta}{2}$. Since
\begin{equation*}
	\begin{aligned} \int_{\mathbb{R}^2}f(x)d\mu_{\sigma}(x)&=\int_{\mathbb{R}^2}f(x)d\mu_{\sigma,n}*\mu_{\sigma,>n}(x)\\ &=\int_{\mathbb{R}^2\times\mathbb{R}^2}f(x+y)d\mu_{\sigma,n}\times \mu_{\sigma,>n}(x,y)\\ &=\int_{\mathbb{R}^2}\int_{\mathbb{R}^2}f(x+y)d\mu_{\sigma,>n}(y)d\mu_{\sigma,n}(x)
	\end{aligned}
\end{equation*}
holds for any $\sigma\in\Omega$ and $n\geq 1$, we obtain
\begin{eqnarray*}
&&\left|\int_{\mathbb{R}^2}f(x)d\mu_{\sigma}(x)-\int_{\mathbb{R}^2}f(x)d\mu_{\sigma,n_0}(x)\right|\nonumber\\
&=&\left|\int_{\mathbb{R}^2}\int_{\mathbb{R}^2}\left(f(x+y)-f(x)\right)d\mu_{\sigma,>n_0}(y)d\mu_{\sigma,n_0}(x)\right|\nonumber\\
&\leq&\int_{B(\textbf{0},\gamma)}\int_{B(\textbf{0},\frac{\gamma}{2^{n_0}})}|f(x+y)-f(x)|d\mu_{\sigma,>n_0}(y)d\mu_{\sigma,n_0}(x) \nonumber \\
&<&\varepsilon.
\end{eqnarray*}

Because $\{\sigma(j)\}_{j=1}^{\infty}$ converges to $\tau$, for any $n\geq1$, there exists $j_0\in\mathbb{N}$ such that for all $j\geq j_0$, we have $\sigma_1(j)\sigma_2(j)\cdots\sigma_n(j)=\tau_1\tau_2\cdots\tau_n$. Hence  $\mu_{\sigma(j),n}=\mu_{\tau,n}$. Therefore, it induces that, for $j\geq j_0$,
\begin{align*}
\left|\int_{\mathbb{R}^2}fd\mu_{\sigma(j)}-\int_{\mathbb{R}^2}fd\mu_{\tau}\right|
&\leq \left|\int_{\mathbb{R}^2}fd\mu_{\sigma(j)}-\int_{\mathbb{R}^2}fd\mu_{\sigma(j),n_0}\right|+\left|\int_{\mathbb{R}^2}fd\mu_{\tau}-\int_{\mathbb{R}^2}fd\mu_{\tau,n_0}\right| \\
&< 2\varepsilon.
\end{align*}
Consequently, $\{\mu_{\sigma(j)}\}_{j=1}^{\infty}$ converges weakly to $\mu_{\tau}$.
\end{proof}

Next, we will consider the situation when the integral periodic zero set of $\mu_{\sigma}$ is an empty set.

\begin{lemma}\label{lem5.4}
Let $\mu_{\sigma}$ be defined by \eqref{eq5.1.1}.  Then  the integral periodic zero set $\textbf{Z}(\mu_\sigma)\neq \emptyset$ if and only if $t_{\sigma_{i}}\equiv t\neq1$ for all $i\geq 1$.
\end{lemma}
\begin{proof}
We first prove the sufficiency. Suppose $t_{\sigma_{i}}\equiv t\neq1$ for all $i\geq 1$, then
$$
\mu_{\sigma}=\delta_{M_1^{-1}t\mathcal{D}}*\delta_{M_1^{-1}M_2^{-1}t
\mathcal{D}}*\cdots*\delta_{M_1^{-1}M_2^{-1}\cdots M_n^{-1}t\mathcal{D}}*\cdots.
$$
Since $M_n\in GL(2, 2\mathbb{Z})$ and $|\det(M_n)|=4$, it follows that $M_n^{*}=2\bar{M}_n$ with  $\bar{M}_n\in GL(2, \mathbb{Z})$ and $|\det(\bar{M}_n)|=1$. Obviously, ${\bar{M}_n}\mathbb{Z}^2=\mathbb{Z}^2$ and ${\bar{M}_n} (2\mathbb{Z}^2)=2{\bar{M}_n} \mathbb{Z}^2=2\mathbb{Z}^2$ for any $n\geq1$.
Then by \eqref{eq2.4}, one may obtain
\begin{equation}\label{eq(3.1)}
	\begin{split}
\mathcal{Z}(\widehat{\mu}_{\sigma})
&=\bigcup_{j=1}^{\infty}M_1^{*}M_2^{*}\cdots M_j^{*}\frac{1}{2t}(\mathring{\mathcal{F}}_2+2\mathbb{Z}^2) \\
&=\frac{1}{t}\bigcup_{j=1}^{\infty}2^{j-1}
\bar{M}_1\bar{M}_2\cdots \bar{M}_j(\mathbb{Z}^2\setminus2\mathbb{Z}^2) \\
&=\frac{1}{t}\bigcup_{j=1}^{\infty}2^{j-1}(\mathbb{Z}^2\setminus2\mathbb{Z}^2) \\
&=\frac{1}{t}(\mathbb{Z}^2\setminus\{\mathbf{0}\}).
\end{split}
\end{equation}
For any $\mathbf{0}\neq\alpha\in \mathcal {F}_{|t|}$, by \eqref{eq(3.1)}, we have $\widehat{\mu}(\frac{\alpha}{t}+k)=0$  for all  $k\in\mathbb{Z}^2$. This shows
$\frac{1}{t}\mathring{\mathcal{F}}_{|t|}\subset\mathbf{Z}(\mu_\sigma)\neq\emptyset$.

Now we prove the necessity.
Suppose that $\textbf{Z}(\mu_\sigma)\neq\emptyset$, then there exists $\xi_0\in\textbf{Z}(\mu_\sigma)$ such that $\widehat{\mu}_\sigma(\xi_0+k)=0$ for all $k\in\mathbb{Z}^2$. It is clear that $\xi_0\notin \mathbb{Z}^2$ since $\widehat{\mu}_\sigma(\xi_0-\xi_0)=\widehat{\mu}_\sigma(\textbf{0})=1\neq 0$.
By \eqref{eq2.4} and \eqref{eq(3.1)}, we have
\begin{equation}\label{eq4.8.0}
\xi_0+\mathbb{Z}^2\subset\bigcup_{j=1}^{\infty}M_1^{*}M_2^{*}\cdots M_j^{*}\frac{1}{2t_{\sigma_j}}
(\mathring{\mathcal{F}}_2+2\mathbb{Z}^2)
=\bigcup_{j=1}^{\infty}\frac{2^{j-1}}{t_{\sigma_j}}(\mathbb{Z}^2\setminus2\mathbb{Z}^2).
\end{equation}

\textbf{Claim:} There exists $1\neq t^*\in \{t_1,t_2,\dots,t_m\}$ such that $\xi_0+\mathbb{Z}^2\subset\frac{\mathbb{Z}^2}{t^*}$.

Obviously, $t^*\neq 1$ since $\xi_0\notin \mathbb{Z}^2$. Noticing $t_{\sigma_j}\in\{t_1,t_2,\dots,t_m\}$, by \eqref{eq4.8.0}, we have
$$\xi_0+\mathbb{Z}^2\subset\bigcup_{j=1}^{m}\frac{\mathbb{Z}^2}{t_j}.$$

We prove the \textbf{Claim} by contradiction.  Suppose there exist $k_i,k'_i,k_j,k'_j\in\mathbb{Z}^2$ such that
$$\xi_0+k_i=\frac{k'_i}{t_i}\quad\text{and}\quad \xi_0+k_j=\frac{k'_j}{t_j}$$
for some $t_i\neq t_j$,
then $\frac{t_j}{t_{i}}k'_i=k'_j+t_j(k_i-k_j)\in\mathbb{Z}^2$. That contradicts $\gcd(t_i,t_{j})=1$ and $\frac{k'_i}{t_{i}}\notin\mathbb{Z}^2$ (because $\xi_0\notin \mathbb{Z}^2$). Hence $\xi_0+\mathbb{Z}^2\subset\frac{\mathbb{Z}^2}{t^*}$ for some $1\neq t^*\in \{t_1,t_2,\dots,t_m\}$.

Let $A:=\{j:t_{\sigma_j}=t^*\}$. If there exists $j_0\geq 1$ such that $t_{\sigma_{j_0}}\neq t^*$, it follows from \textbf{Claim} and  \eqref{eq4.8.0}   that
\begin{equation*}
\begin{split}
	\xi_0+\mathbb{Z}^2
	&\subset \bigcup_{j\in A}\frac{2^{j-1}}{t_{\sigma_j}}(\mathbb{Z}^2\setminus2\mathbb{Z}^2) \\
	&\subset \frac{1}{t^*}\bigcup_{j=1,j\neq j_0}^{\infty}2^{j-1}(\mathbb{Z}^2\setminus2\mathbb{Z}^2) \\
	&\subset \frac{1}{t^*}(\mathbb{Z}^2\setminus 2^{j_{0}-1}(\mathbb{Z}^2\setminus2\mathbb{Z}^2)).
\end{split}
\end{equation*}

Next, we will show that there exists $k_0\in \mathbb{Z}^2$ such that
\begin{equation}\label{equ-4.5}
	 \xi_0+k_0\in \frac{1}{t^*}(  2^{j_{0}-1}(\mathbb{Z}^2\setminus2\mathbb{Z}^2)).
\end{equation}
Noting that $t^*\in 2\mathbb{Z}+1$ and $\xi_0\in \frac{\mathbb{Z}^2}{t^*}$, let $t^*=2^{r_0}l+1$ and $\xi_0:=2^{s_0}\frac{v_0}{t^*}$, where $r_0>0$, $l\in2\mathbb{Z}+1$, $s_0\geq0$  and $ v_0\in \mathbb{Z}^2\setminus2\mathbb{Z}^2$. 
We consider the following two cases:

(i) If $s_0=0$, let $m$ be a sufficiently large
positive integer such that ${t^*}^{2^m}=(2^{r_0}l+1)^{2^m}=2^NL+1$ with $N>j_0$, and let $t'=2^{j_{0}-1}-1$.  It is easy to get that ${t^*}^{2^m}t'=2^{j_{0}-1}\theta-1$ for some $\theta\in 2\mathbb{Z}+1$. Let $k_0={t^*}^{2^m-1}t'v_0$, hence
$$
\xi_0+k_0=\frac{v_0}{t^*}+\frac{{t^*}^{2^m}t'v_0}{t^*}=\frac{2^{j_{0}-1}\theta v_0}{t^*}\in \frac{1}{t^*}(  2^{j_{0}-1}(\mathbb{Z}^2\setminus2\mathbb{Z}^2)).
$$
Hence \eqref{equ-4.5} holds.

(ii) If $s_0\neq 0$, let $k_0=({t^*}^{2^m-1}t'(2^{s_0}+t^*)+1)v_0$.  Similar to the case (i), we  have
$$
\xi_0+k_0=\frac{2^{j_{0}-1}(2^{s_0}+t^*)\theta v_0}{t^*}\in \frac{1}{t^*}(  2^{j_{0}-1}(\mathbb{Z}^2\setminus2\mathbb{Z}^2)),
$$ 

yielding \eqref{equ-4.5}.

Combining (i) and (ii), it concludes that $t_{\sigma_{i}}\equiv t\neq1$ for all $i\geq 1$.
So we complete the proof.
\end{proof}

At the final part of this section, we first prove Theorem \ref{th1.6}, then prove Theorem \ref{th1.5}.

\begin{proof}[\textbf{Proof of Theorem \ref{th1.6}}]
For the sufficiency. Suppose that $t_2\mid t_1$, then we take $r=\frac{t_1}{t_2}\in 2\mathbb{Z}+1$ and let $\widetilde{D}_n=\frac{1}{t_2}D_n$ for $n\geq 1$. By Lemma \ref{lem2.2}, we know that $\mu_{\{M_n\},\{D_n\}}$ is a spectral measure if and only if
\begin{equation*}
\mu_{\{M_n\},\{\widetilde{D}_n\}}=\delta_{M_1^{-1}r\mathcal{D}}*\delta_{M_1^{-1}M_2^{-1}\mathcal{D}}*\delta_{M_1^{-1}M_2^{-2}\mathcal{D}}*\cdots=\delta_{M_1^{-1}\widetilde{D}_1}*v_1\circ M_1
\end{equation*}
is a spectral measure, where $v_1$ is defined by \eqref{eq1.2.0}.

Since $M_2\in GL(2,2\mathbb{Z})$, the measure $v_1$  is a spectral measure by Theorem \ref{th1.2.9}. We claim that $\mathbb{Z}^2$ is the unique spectrum of $v_1$ containing $\mathbf{0}$.
By \eqref{eq(3.1)}, it is easy to see that $\mathcal{Z}(\widehat{v}_1)=\mathbb{Z}^2\setminus\{\mathbf{0}\}$.
Let $\Gamma$ be a spectrum of $v_1$ with $\mathbf{0}\in\Gamma$, then we know $\Gamma\subset\mathcal{Z}(\widehat{v}_1)\cup\{\mathbf{0}\}=\mathbb{Z}^2$ by \eqref{eq2.5}. If $\Gamma\neq\mathbb{Z}^2$, then  $a\bot\Gamma$ holds for any $a\in \mathbb{Z}^2\setminus\Gamma$. This is a contradiction, thus the claim follows.

Let $L=\frac{1}{2}M_1^{*}\{(0,0)^t,(1,0)^t,(0,1)^t,(1,1)^t\}$, then $(M_1,r\mathcal{D},L)$ is a Hadamard triple.  Now we construct a set $\Lambda=L+M_1^{*}\mathbb{Z}^2$. For any $\xi\in\mathbb{R}^2$, we apply Theorem \ref{thm-JP98} (ii) to get
\begin{eqnarray*}
Q_{\mu_{\{M_n\},\{\widetilde{D}_n\}},\Lambda}(\xi)&=&\sum\limits_{\lambda\in\Lambda}\left|\widehat{\mu}_{\{M_n\},\{\widetilde{D}_n\}}(\xi+\lambda)\right|^2\\
&=&\sum\limits_{l\in L}\sum\limits_{k\in\mathbb{Z}^2}\left|\widehat{\mu}_{\{M_n\},\{\widetilde{D}_n\}}(\xi+l+M_1^{*}k)\right|^2\\
&=&\sum\limits_{l\in L}\sum\limits_{k\in\mathbb{Z}^2}\left|m_{r\mathcal{D}}(M_1^{*-1}(\xi+l)+k)\right|^2\left|\widehat{v}_1(M_1^{*-1}(\xi+l)+k)\right|^2\\
&=&\sum\limits_{l\in L}\left|m_{r\mathcal{D}}(M_1^{*-1}(\xi+l))\right|^2\sum\limits_{k\in\mathbb{Z}^2}\left|\widehat{v}_1(M_1^{*-1}(\xi+l)+k)\right|^2\\
&=&\sum\limits_{l\in L}\left
|m_{r\mathcal{D}}(M_1^{*-1}(\xi+l))\right|^2=1.
\end{eqnarray*}
Hence $\mu_{\{M_n\},\{\tilde{D}_n\}}$ is a spectral measure. The sufficiency is proved.

For the necessity. Assume that $\mu_{\{M_n\},\{D_n\}}$ is a spectral measure. Analogous to \eqref{eq3.1},
$\mu_{\{M_n\},\{D_n\}}=\mu_{\{\widetilde{M}_n\},\mathcal {D}}$, where $\widetilde{M}_1=\frac{1}{t_1}M_1$, $\widetilde{M}_2=\frac{t_1}{t_2}M_2$ and  $\widetilde{M}_n=M_2, n\geq3$.  Let $p=\textup{lcm}(t_1,t_2)$ and $\Lambda$ containing $\mathbf{0}$ be a spectrum of $\mu_{\{\widetilde{M}_n\},\mathcal {D}}$. For any $\mathbf{0}\neq s\in \mathcal{F}_p$, choose $l_s\in\mathcal{F}_2$ such that $s+pl_s\notin 2\mathbb{Z}^2$ and $l_\textbf{0}=\textbf{0}$. It follows from Lemma \ref{lem4.2} that
\begin{equation*}
\Gamma=\bigcup\limits_{s\in \mathcal{F}_p}\left(\frac{s+pl_s}{2p}+\Lambda_{s,l_s}\right)
\end{equation*}
is a spectrum of $\mu_1$,  where $\Lambda_{s,l_s}$ and
$\mu_1$ are defined by \eqref{eq-3.6.0} and \eqref{eq4.8.1}, respectively.
Noting that $M_2\in GL(2,2\mathbb{Z})$ and $\det(M_2)=4$, similar to \eqref{eq(3.1)},
by \eqref{eq2.4}, we have
\begin{equation*}
\mathcal{Z}(\widehat{\mu}_1)=\bigcup_{j=2}^{\infty}\widetilde{M}_{2}^{*}\cdots \widetilde{M}_{j}^{*}\mathcal{Z}(m_{\mathcal {D}})=\frac{t_1}{t_2}\bigcup_{j=1}^{\infty}M_2^{*j}\mathcal{Z}(m_{\mathcal {D}})
=\frac{t_1}{t_2}(\mathbb{Z}^2\setminus\{\textbf{0}\}).
\end{equation*}

It is easy to show that $\frac{t_1}{t_2}\mathbb{Z}^2$ is the unique spectrum of $\mu_1$ containing $\textbf{0}$. This implies that
$$
\frac{t_1}{t_2}\mathbb{Z}^2=\bigcup\limits_{s\in \mathcal{F}_p}\left(\frac{s+pl_s}{2p}+\Lambda_{s,l_s}\right).
$$
Combining $s+pl_s\notin 2\mathbb{Z}^2$ for $s\neq\textbf{0}$ and $t_1, t_2 \in 2\mathbb{Z}+1$, it yields that $\Lambda_{s,l_s}=\emptyset$ for all
$s\neq \textbf{0}$. Hence
$$
\frac{t_1}{t_2}\mathbb{Z}^2=\Lambda_{\textbf{0},\textbf{0}}\subset\mathbb{Z}^2.
$$
Therefore, $t_2\mid t_1$, and the necessity is proved.
\end{proof}

\begin{cor}\label{cor5.1}
Given an integer  $m\geq2$,  let $\{M_n\}_{n=1}^{m}\subset GL(2, \mathbb{Z})$ be a finitely many expanding matrices with $|\det(M_n)|=4$. For $n\geq m$, let $M_n\equiv M_m$ and let digit sets
\begin{equation*}
D_n=\left\{
\begin{aligned}
&t_n\mathcal{D} \ \ 1\leq n\leq m-1\\
&t_m\mathcal{D} \ \ n\geq m
\end{aligned}\right. \quad (t_n\in2\mathbb{Z}+1, n=1,2,\dots,m.)
\end{equation*}
with $\mathcal{D}=\{(0,0)^t,(1,0)^t,(0,1)^t,(-1,-1)^t\}$. If $\mu_{\{M_n\},\{D_n\}}$ defined by \eqref{eq1.2} is a spectral measure, then $t_m\mid t_{m-1}$.
\end{cor}

\begin{proof}
 If $\mu_{\{M_n\},\{D_n\}}$ is a spectral measure, 	from Theorem \ref{th1.1}, we have that $M_n\in GL(2,2\Z)$ for all $n\ge 2$. 	Therefore, together with  Remark \ref{re3.5}, the proof  of the necessity of Theorem \ref{th1.6} also yields the corollary directly.
\end{proof}

\bigskip

Finally, we give the proof of Theorem \ref{th1.5} which is inspired by Wu and Xiao \cite{Wu-Xiao2022}. For $\sigma=(\sigma_n)_{n=1}^{\infty}\in\Omega$, we define the shift operator
$$\varrho^{n}(\sigma)=\sigma_{n+1}\sigma_{n+2}\sigma_{n+3}\cdots$$
and let
$$
v_ {\sigma,n}=\delta_{M_{{n+1}}^{-1}D_{\sigma_{n+1}}}*
\delta_{M_{{n+1}}^{-1}M_{{n+2}}^{-1}D_{\sigma_{n+2}}}*\cdots
$$
for convenience.

\bigskip

\begin{proof}[\textbf{Proof of Theorem \ref{th1.5}}]
We first prove  the sufficiency, i.e., $\mu_{\sigma}$ is a spectral measure if $\sigma\notin \cup_{l=1}^{\infty}\Sigma_l$, where $\Sigma_l=\{i_1i_2\cdots i_lj^{\infty}\in\Omega:i_l\neq j, j\neq1\}$.  To achieve this goal, it suffices, by Theorem \ref{th2.5}, to identify a subsequence of $\{\varrho^n(\sigma)\}$ that converges to some $\tau$ such that $\mathbf{Z}(\mu_{\tau}) = \emptyset$. We divide the proof into the following two cases.

\textbf{Case i:}  `$1$' appears infinitely many times in  $\sigma=(\sigma_n)_{n=1}^\infty$.

Let $\{n_k\}_{k=1}^\infty$ be a  subsequence  of $\{n\}_{n=1}^\infty$ such that $\sigma_{n_k}=1$ for all $k$.
It is clear that $\{\varrho^{n_k-1}(\sigma)\}_{k=1}^{\infty}\subset\Omega$, and for all $k\geq1$ there exists $\tau_k\in\Omega$ such that $\varrho^{n_k-1}(\sigma)=1\tau_k$. Due to the compactness of $\Omega$, there must exist a subsequence $\{k_j\}_{j=1}^{\infty}$ and $\widetilde{\tau}\in\Omega$ such that $\{\varrho^{n_{k_j}-1}(\sigma)\}_{j=1}^{\infty}$ converges to $\tau=1\widetilde{\tau}$. Since $\{v_ {\sigma,n_{k_j}-1}\}=\{\mu_{\varrho^{n_{k_j}-1}(\sigma)}\}$, we know that $v_{\sigma,n_{k_j}-1}$ converges weakly to $\mu_{\tau}$ and $\textbf{Z}(\mu_{\tau})=\emptyset$ by Lemmas \ref{lem5.2.0} and  \ref{lem5.4}. Thus $\mu_{\sigma}$ is a spectral measure by  Theorem \ref{th2.5}.

\textbf{Case ii:} `$1$' appears finitely many times in $\sigma=(\sigma_n)_{n=1}^\infty$.

From $\sigma\notin \cup_{l=1}^{\infty}\Sigma_l$, it follows that there must exist $j_1\neq j_2\in \{2,3,\dots,m\}$ such that $j_1$ and $j_2$ appear infinitely many times in  $\sigma=(\sigma_n)_{n=1}^\infty$. In ascending order, we let
$$\{l_1^{(i)},l_2^{(i)},l_3^{(i)},\dots\}=\{n\geq1:\sigma_n=j_i\}\quad\text{where }  i=1,2.$$
According to the gaps between $l_{k}^{(i)}$ and $l_{k+1}^{(i)}$, this case can be considered from two aspects.

On the one hand, if $\varlimsup\limits_{k\rightarrow\infty}\left(l_{k+1}^{(i)}-l_{k}^{(i)}\right)<\infty$ for all $i=1,2$. Since $\Omega$ is a compact space and  $\{\varrho^{k}(\sigma)\}_{k=1}^{\infty}\subset\Omega$, there exist a subsequence $\{k_j\}_{j=1}^{\infty}$ and $\eta\in\Omega$ such that $\{\varrho^ {k_j}(\sigma)\}_{j=1}^{\infty}$ converges to $\eta$ and $j_1,j_2$ also appear infinitely many times in $\eta$.  Lemmas \ref{lem5.2.0} and \ref{lem5.4} yield that $v_{\sigma,k_j}$ converges weakly to $\mu_{\eta}$ with $\textbf{Z}(\mu_{\eta})=\emptyset$. Therefore, $\mu_{\sigma}$ is a spectral measure by Theorem \ref{th2.5}.

On the other hand, if $\varlimsup\limits_{k\rightarrow\infty}\left(l_{k+1}^{(i)}-l_{k}^{(i)}\right)=\infty$ for some $i\in\{1, 2\}$. Without loss of generality, we assume that $i=1$.  Then there   exists a subsequence $\{k_s\}_{s=1}^{\infty}$ such that
$\lim\limits_{s\rightarrow\infty}\left(l_{k_{s}+1}^{(1)}-l_{k_{s}}^{(1)}\right)
=\infty$.  In this situation, there is some $\zeta_s\in\Omega$ such that $\varrho^{l_{k_s}^{(1)}-1}(\sigma)=j_1\zeta_s$ for all $s\geq1$. Similarly,   there exists a subsequence $\{\widetilde{k}_s\}_{s=1}^{\infty}\subset\{k_s\}_{s=1}^{\infty}$ and $\omega=j_1\zeta,\zeta\in\{1,2,\dots,j_1-1,j_1+1,\dots,m\}^{\mathbb{N}}$ such that $\varrho^{l_{\widetilde{k}_s}^{(1)}-1}(\sigma)$ converges to $\omega$. So  $v_{\sigma,l_{\widetilde{k}_s}^{(1)}-1}$ converges weakly to $\mu_{\omega}$ with $\textbf{Z}(\mu_{\omega})=\emptyset$ by Lemmas \ref{lem5.2.0} and \ref{lem5.4}. Therefore, the desired result follows by Theorem \ref{th2.5}.

Finally, for the necessity, it suffices to show that if $\sigma\in \cup_{l=1}^{\infty}\Sigma_l$, then $\mu_{\sigma}$ is not a spectral measure. Under the assumption, there exists some $l_1\geq1$ such that $\sigma\in \Sigma_{l_1}=\{i_1i_2\cdots i_{l_1}j^{\infty}: i_{l_1}\ne j, j\ne 1\}$.  Since $\gcd(t_{\sigma_{l_1}},t_{\sigma_{l_1+1}})=1$ and $t_{\sigma_{l_1+1}}\neq1$, it concludes that $\mu_{\sigma}$ is not a spectral measure by Corollary \ref{cor5.1}.

The proof of Theorem \ref{th1.5} is finished.
\end{proof}

\bigskip

\section{Concluding remarks}\label{sect.5}
In the final section, we present some remarks and pose open questions that merit further investigation. Notably, Theorem \ref{th1.1} in this paper addresses only the case where $\alpha_{n_1}\beta_{n_2}-\alpha_{n_2}\beta_{n_1} \notin 2\mathbb{Z}$. However, it is potentially more intriguing to explore the spectrality of $\mu_{\{M_n\},\{D_n\}}$ when this determinant condition lies within $2\mathbb{Z}$.

For instance, consider the setting where $M_{2k-1} = 4I$, $M_{2k} = 3I$, $D_{2k-1} = 2\mathcal{D}$, and $D_{2k} = \mathcal{D}$ for all $k \geq 1$, with
$$
\mathcal{D} = \left\{(0,0)^t, (1,0)^t, (0,1)^t, (-1,-1)^t \right\}.
$$
In this case, the Moran measure $\mu_{\{M_n\},\{D_n\}}$ takes the form:
$$
\mu_{\{M_n\},\{D_n\}} = \delta_{\frac{\mathcal{D}}{2}} * \delta_{\frac{\mathcal{D}}{12}} * \delta_{\frac{\mathcal{D}}{24}} * \delta_{\frac{\mathcal{D}}{144}} * \cdots,
$$
which can be regrouped into the self-affine measure
$$
\mu_{M,D} := \delta_{\frac{\mathcal{D} + 6\mathcal{D}}{12}} * \delta_{\frac{\mathcal{D} + 6\mathcal{D}}{144}} * \cdots
$$
with $M = 12I$ and $D = \mathcal{D} + 6\mathcal{D}$. It is readily verified that
$$
\frac{1}{4}\mathring{\mathcal{F}}_4 \subset \mathcal{Z}(m_D) = \left( \frac{1}{2} \mathring{\mathcal{F}}_2 + \mathbb{Z}^2 \right) \cup \frac{1}{6} \left( \frac{1}{2} \mathring{\mathcal{F}}_2 + \mathbb{Z}^2 \right),
$$
where $\mathring{\mathcal{F}}_n$ is defined in \eqref{eq-3-2-1}. Consequently, the set $L := \frac{1}{4}M^* \mathcal{F}_4 \subset \mathbb{Z}^2$, and the triple $(M, D, L)$ forms a Hadamard triple. By a result of Dutkay \emph{et al.} \cite{Dutkay-Haussermann-Lai2019}, this implies that $\mu_{M,D} = \mu_{\{M_n\},\{D_n\}}$ is a spectral measure.

More generally, consider the case where
$$
M_{2k-1} = \begin{pmatrix}
	4 & 0 \\
	0 & 6
\end{pmatrix}, \
M_{2k} = \begin{pmatrix}
	3 & 0 \\
	0 & 2
\end{pmatrix}, \
D_{2k-1} = \begin{pmatrix}
	2 & 0 \\
	0 & 3
\end{pmatrix} \mathcal{D}, \
D_{2k} = \mathcal{D}, \quad \text{for } k \ge 1,
$$
then the associated Moran measure satisfies $\mu_{\{M_n\},\{D_n\}} = \mu_{12I, \mathcal{D} + 6\mathcal{D}}$, which is again a spectral measure. This observation motivates the following natural question:
\begin{qu}
	For the Moran measure $\mu_{\{M_n\},\{D_n\}}$ considered in Theorem \ref{th1.1}, if the condition $\alpha_{n_1} \beta_{n_2} - \alpha_{n_2} \beta_{n_1} \notin 2\mathbb{Z}$ is removed, can one establish a necessary and sufficient condition under which $\mu_{\{M_n\},\{D_n\}}$ is a spectral measure?
\end{qu}

In \cite{Deng-Li2022}, Deng and Li investigated the Moran measure $\mu_{\{p_n\},\{D_n\}}$ on $\mathbb{R}$, where $D_n = \{0, d_n\}$ and $p_n, d_n$ are integers satisfying $|p_n| > |d_n| > 0$, with $\{d_n\}$ bounded. They proved that $\mu_{\{p_n\},\{D_n\}}$  is a spectral measure if and only if the number of factors of $2$ in the sequence $\left\{\frac{p_1p_2\cdots p_n}{2d_n}\right\}$ are distinct. However, when extended to $\mathbb{R}^2$, the situation becomes significantly more complex and challenging to analyze.

More recently, Wu \cite{Wu2024} established a necessary and sufficient condition for the spectrality of self-similar measures on $\mathbb{R}$ generated by alternating contraction ratios. Luo \emph{et al.} \cite{Luo-Mao-Liu2024} further studied a class of Moran measures with alternating contraction ratios, extending the main result of \cite{An-He2014}. Inspired by these developments, it would be a worthwhile direction to explore similar extensions to higher-dimensional settings.

In Theorem \ref{th1.5}, we assumed that $|\det(M_n)| = 4$ and $D_n = t_n\mathcal{D}$ for all $n \ge 1$. Under these assumptions, the associated Moran set
$$
T = \sum_{n=1}^{\infty} M_1^{-1} M_2^{-1} \cdots M_n^{-1} D_n
$$
may have positive Lebesgue measure. In this context, it is natural to consider the relevance of Fuglede's spectral set conjecture to the set $T$; see, for example, \cite{Liu-Liu-Tang2023}.

Moreover, for general digit sets $\{D_n\}$ of the form given in \eqref{eq1.3}, it remains an open and meaningful problem to study the spectral properties of $\mu_{\{M_n\},\{D_n\}}$ when $|\det(M_n)| = 4$ for all $n \ge 1$.

\end{document}